# Integrated Neighborhood Colorings of Graphs

Robert Cowen
Department of Mathematics[1]
Queens College, CUNY

Abstract. . The idea that those that are different from you are "unfriendly" is captured in the definition of unfriendly 2-colorings in graph theory [ 1 ] and in the definition of unfriendly n-partitions [ 9 ]. It has been proved that every finite graph has an unfriendly 2-coloring [ 1 ]. We give a more general definition for k ⩾ 2 that we call "integrated" rather than "unfriendly." Then we prove that every finite graph has an integrated k-coloring. The two definitions coincide, however, when k = 2. We then give some applications to general graph coloring problems as well as max-cut problems.

Introduction. Unfriendly colorings of graphs were introduced about thirty years ago and it was proved that every finite graph has an unfriendly 2-coloring ([ 1 ]). Somewhat stronger results were proved even earlier by Bernardi [ 2 ] and others (see the references in [ 2 ]). In [ 5 ] we showed how to find unfriendly 2-colorings in polynomial time with a computer using *Mathematica*. In [ 9 ], a generalization to k-colorings, called "unfriendly k-partitions" was introduced primarily to study infinite graphs. Here we give a generalization to k-colorings that differs from the one given in [ 9 ]. Our definition, called "integrated neighborhood k-colorings," leads to new applications, when k > 2. As in the k = 2 case, computer programs can be written to find these integrated colorings in polynomial time. Finally we prove a compactness result extending our main theorems to locally finite infinite graphs.

Preliminaries. A *k-coloring* of a graph $G = <V, E>$ is a function $c: V \to \{1, \ldots, k\}$. A k-coloring is called *proper* if $c(x) \neq c(y)$, if x and y are adjacent (connected by an edge). A *neighborhood of a vertex v, N(v)* consists of all vertices that adjacent to *v*. Note that, according to this definition, *v* is not a member of *N*(*v*). Suppose that graph *G* has a k-coloring. Then *N*(*v*) will be called an *integrated neighborhood* of the colored graph *G*, if at most $(1/k)|N(v)|$ vertices in *N*(*v*) have the same color as *v*, where $|N(v)|$ is the cardinality of *N*(*v*). Thus, if *N*(*v*) is not an integrated neighborhood, then more than $(1/k)|N(v)|$ of its vertices have the same color as *v*. If there is a k-coloring such that for every vertex *v*, N(*v*) is an integrated neighborhood, the coloring will be called an *integrated neighborhood k-coloring* of the graph *G* or just simply an *integrated k-coloring*. We note that our definition of an integrated k-coloring differs from the definition of an unfriendly k-partition of a graph $G = <V, E>$ as defined in [ 9 ]; $c: V \to \{1, \ldots, k\}$, is an *unfriendly k-partition* if for every vertex v, $|\ y \in N(v): c(y) = c(v)| \leq |y \in N(v): c(y) \neq c(v)|$, that is, at most $(1/2)|N(v)|$ vertices in *N*(*v*) have the same color as *v*. However, If k = 2, the two definitions coincide. If *G* is a finite graph, $\Delta(G)$ is the maximum degree of the vertices in *G*. An edge of a colored graph is called <u>mixed</u>, if its endpoints have different colors. The <u>mixing number</u> of a colored graph is the number of its mixed edges.

Main Results. We shall prove first that every finite graph *G* has an integrated k-coloring, 2 ≤ k ≤ n, where n is the number of vertices of *G*. We first prove a more general result from an unpublished paper of ours with William Emerson [ 6 ].

Let $c: V \to \{1, ..., k\}$ be a k-coloring of the graph $G = <V, E>$ and let $d_i(v)$ denote the number of neighbors of vertex v with color i, and let $d(v)$ denote the degree of v. The following result is from our unpublished paper with William Emerson with some improvements suggested by Douglas Woodall.

<u>Lemma</u>. Let $G = <V,E>$ and let $p_1, ..., p_k$ be such that $0 \leq p_i \leq 1$, $1 \leq i \leq k$ and $\sum_{i=1}^{k} p_i \geq 1$. Then there is a k-coloring $c : V \to \{1, ..., k\}$ such that $d_{c(v)}(v) \leq p_{c(v)} d(v)$, for all v in V.

Proof. Let c be a k-coloring of G that minimizes $\sigma_c = \sum_{i=1}^{k} b_i / p_i$, where $b_i$ is the number of edges of G with both endpoints colored i. If there was a vertex $v_0$ with $c(v_0) = i_0$, such that $d_{i_0}(v_0) > p_{i_0} d(v_0)$, then, we claim there must exists a $j_0$, $j_0 \neq i_0$, $1 \leq j_0 \leq k$, such that, $d_{j_0}(v_0) < p_{j_0} d(v_0)$; since if $d_i(v_0) \geq p_i d(v_0)$, $i \neq i_0$, then,

$$d(v_0) = d_{i_0}(v_0) + \sum_{i \neq i_0} d_i(v_0) > p_{i_0} d(v_0) + \sum_{i \neq i_0} p_i d(v_0) = d(v_0) \sum_{i=1}^{k} p_i = d(v_0)!$$

Thus, there exists a j, $j \neq i_0$, $1 \leq j_0 \leq k$, such that, $d_{j_0}(v_0) < p_{j_0} d(v_0)$. So $d_{i_0}(v_0) > p_{i_0} d(v_0)$ and $d_{j_0}(v_0) < p_{j_0} d(v_0)$. Therefore, $(d_{i_0}(v_0)/ p_{i_0}) > (d_{j_0}(v_0)/p_{j_0})$. Hence, if we changed $v_0$'s color from $i_0$ to $j_0$ it would decrease $\sigma_c$ which is impossible, by the minimality of $\sigma_c$.

<u>Theorem 1</u>. Let G be a finite graph with n vertices and suppose $2 \leq k \leq n$. Then G has an integrated k-coloring.

Proof. Let $p_i = 1/k$, $1 \leq i \leq k$. Then $\sum_{i=1}^{k} p_i = 1$. So, by the lemma, there is a k-coloring $c : V \to \{1, ..., k\}$ such that $d_{c(v)}(v) \leq (1/k)d(v)$ for all v in V.

<u>Corollary 1.1</u>. Let G be a finite graph with n vertices and suppose $2 \leq k \leq n$. Then G has a k-coloring such that at most $(1/k)|N(v)|$ vertices adjacent to v have the same color as v.

<u>Finding an Integrated k-Coloring in Polynomial Time</u>. In [5] we gave a computer program for finding an unfriendly 2-coloring for a finite graph. We show here that essentially the same ideas can be used to give a polynomial time algorithm for finding an integrated k-coloring of a finite graph. Call a vertex *k-secure* if more than 1/k of its neighbors have the same color as the vertex. We claim that the following will produce an integrated k-coloring in polynomial time, O(ED), where E is the number of edges and D is the maximum degree of G.

1. Find a k-secure vertex, if one exists; if there is none, we have an integrated k-coloring. STOP
2. If a k-secure vertex v is found under step 1, there must be some color $i_0$, $1 \leq i_0 \leq k$, with less than $(1/k) |N(v)|$ vertices in N(v) with the color $i_0$. Change the color of v to $i_0$. Return to step 1.

In step 2, if a k-secure vertex v is found, it is necessary to count the number of vertices of v's color versus those vertices of other colors; this takes O(D) time where D is the maximum degree of G. The number of times it is necessary to do this is, in the worst case, E, the number of edges in G, since each time we change v's color the mixing number increases and the mixing number is bounded by E. So the algorithm runs in O(ED) time.

Shelah and Milner [9] prove that every graph, finite or infinite has an unfriendly 3-partition; this implies that every finite graph has a 3-coloring such that at most $(1/2) |N(v)|$ vertices in N(v) have the same color as v. This corollary shows that for finite graphs, this can be sharpened to at most $(1/3) |N(v)|$ vertices in N(v) have the same color as v.

Suppose now that that $k = \Delta(G) + 1$, where $\Delta(G)$ is the maximum degree of the vertices in G. Then $|N(v)| < k$, for all v ∈ V. In this case $(1/k)|N(v)| < 1$, that is, the k-coloring is a proper coloring and we have a proof of the following result of Brooks [ 3 ].

<u>Corollary 1.2</u>. Let G be a finite graph and let $k = \Delta(G) + 1$, where $\Delta(G)$ is the maximum degree of the vertices in G. Then $\chi(G) \leq \Delta(G) + 1$, where $\chi(G)$ is the chromatic number of G.

In the mid 1980s several authors "relaxed" the requirement of proper colorings that vertices that share an edge must be colored differently. One such "relaxation" was called "defective colorings." According to [4], a (k,u)-coloring of a graph G is a k-coloring of its vertices such that for any vertex v, at most u of its neighbors have the same color as v. (An (k, 0)- coloring is the same as a proper k-coloring.) The following result from [ 7 ], now follows easily from Corollary 1.1; in addition we point out that a (k,u)-coloring can be found in polynomial time.

<u>Corollary 1.3</u>. Let G be a finite graph and let $u \geq 0$ and $k \geq 1$ be integers such that $\Delta(G) \leq k(u + 1) - 1$. Then G is (k, u) - colorable. Moreover, a (k,u)-coloring can be found in polynomial time.

<u>Proof</u>. Let c be an integrated k-coloring of G. We claim that c is also a (k, u)-coloring of G. Corollary 1.1 implies that at most $(1/k)|N(v)|$ vertices adjacent to v have the same color as v. Clearly $|N(v)| \leq \Delta(G)$, for all vertices v of G. However, $(1/k)|N(v)| \leq (1/k) \Delta(G) \leq (1/k)(k(u + 1) - 1) = (u + (1 - 1/k) < u + 1$. Thus at most u vertices adjacent to v have the same color as v under the coloring c. Moreover, finding an integrated k-coloring for a finite graph can easily be carried out in polynomial time as we explained above.

<u>Theorem 2</u>. Let G be a finite graph with n vertices and m edges and let $2 \leq k \leq n$. In any integrated k-coloring of graph G at least $((k-1)/k)m$ of the edges are mixed.

<u>Proof</u>. Suppose G has a fixed integrated k-coloring. Let M be the mixing number of this coloring. Since the k-coloring is integrated, for each vertex v, at most $(1/k)|N(v)|$ are colored with the same color as v; so, at least $(k - 1)/k|N(v)|$ of the the vertices of N(v) are colored differently than v. Therefore at least $(k - 1)/k|N(v)|$ edges with v as one of its endpoints are mixed; let the exact number of mixed edges containing v be denoted by mix(v). Therefore, $mix(v) \geq ((k - 1)/k)|N(v)|$. Now $\Sigma mix(v) = 2 M$, where the sum is taken over all vertices v, since each mixed edge is counted twice, once at each endpoint. Also summing over all vertices v, we get,

$\quad 2M = \Sigma mix(v) \geq \Sigma ((k - 1)/k)|N(v)| = (k - 1)/k \Sigma |N(v)| = (k - 1)/k (2 m)$,

since $\Sigma |N(v)|$ counts each edge twice. Hence $M \geq ((k - 1)/k) m$, as required.

<u>Corollary 1</u>. Let G be a finite graph with n vertices and m edges and let $2 \leq k \leq n$. Then there is a k-coloring where at least $((k-1)/k)m$ of the edges are mixed.

<u>Proof</u>. Follows from Theorem 1 and Theorem 2.

<u>Corollary 2</u>. Any finite graph can be 2-colored so that at least 1/2 of its edges are mixed.

<u>Proof</u>. Follows from Corollary 1 with k = 2.

As we have already mentioned, this coloring can be found in polynomial time .

The *max-cut problem* is to partition the vertices of a graph *G* into two sets so as to maximize the number of edges with vertices in both sets. This is an NP-hard problem; however, it is easy to show this maximum number is at least m/2 and to find a cut with at least m/2 edges.

Corollary 3. Let G be a finite graph with m edges. Then the vertices of G can be partitioned into 2 sets $V_1$ and $V_2$ with at least m/2 edges with vertices in both sets. Moreover such a cut can be found in polynomial time.

Proof. Follows from Corollary 2 by taking $V_1$ to be the set of vertices receiving one of the colors and $V_2$ to be the set of vertices receiving the other color. The computer program in [5] will give a polynomial time solution.

The *k-max-cut problem* to partition the vertices of a graph *G* into k-sets so as to maximize the number of edges with endpoints in two different sets. We show next that this maximum number is at least ((k-1)/k )m.

Corollary 4. Let G be a finite graph with m edges. The the vertices of G can be partitioned into k sets: $V_1$, ... , $V_k$ with at least ((k-1)/k )m edges with vertices in different sets of the partition.

Proof. This follows from Corollary 1, by taking $V_i$ to be the set of vertices that receive the *i-th* color, $1 \leq i \leq k$.

Finally, we extend our result in Theorem 1 to locally finite infinite graphs, that is, to those infinite graphs in which every vertex has only finitely many neighbors. We make use of the Rado Selection Lemma[8].

Rado Selection Lemma. Let { $A_v$: v∈ I} be a set of finite sets. Suppose for every finite W ⊂ I, there exists a function, $f_W$, with domain W, such that $f_W(v) \in A_v$, v ∈ W. Then there exist a function *f*, with domain I such that for every finite W, there exists a finite $W_1 \supset W$ with $f(v) = f_{W_1}(v)$, v ∈ W.

Theorem 3. Let *G* be a locally finite graph. Suppose k is a positive integer. Then *G* has an integrated k-coloring.

Proof. Let I be the set of vertices of *G* and if v ∈ I, let $A_v$ = {1,2, ... ,k} (the set of "colors"}. For each finite set of vertices, W, there is an integrated coloring, $f_W$, of the subgraph $G_W$, induced in G by vertex set W, by Theorem 1. Then $f_W(v) \in A_v$, v ∈ W. The Rado Selection Lemma gives a function *f* with domain I such that for every finite W, there exists a finite $W_1 \supset W$ with f(v) = $f_{W_1}(v)$, v ∈ W. We claim *f* is an integrated coloring of *G*. Suppose *v* is a vertex of *G* and N(v) is the set of neighbors of *v*. We must show that at most (1/k)|N(v)| vertices of N(v) have the same color as *v* under the coloring *f*. Let W be the set of vertices in N(v); by the Rado Selection Lemma, there exists a finite $W_1 \supset W$ with $f(v) = f_{W_1}(v)$, v ∈ W. Since $f_{W_1}$ is an integrated coloring, at most (1/k)|N(v)| vertices of N(v) have the same color as *v* under the coloring $f_{W_1}$; since $f(v) = f_{W_1}(v)$, v ∈ N(v), the same is true for *f*.

Since Theorem 1 can be extended to locally finite infinite graphs, so can all of its corollaries. In particular, we have the following.

<u>Corollary 3.1</u>. Let *G* be a locally finite graph and let u ⩾ 0 and k ⩾ 1 be integers such that $\Delta(G) \leq k(u+1) - 1$. Then G is (k, u) - colorable.

[1]Professor Emeritus